\documentclass{amsart}
\usepackage{amscd}
\usepackage{amssymb}

\newcommand{\F}{{\mathbb F}}
\newcommand{\Q}{{\mathbb Q}}

\newcommand{\Gal}{\operatorname{Gal}\,}
\newcommand{\GL}{\operatorname{GL}}

\newcommand{\Cl}{\operatorname{Cl}}

\newcommand{\lra}{\longrightarrow}

\newcommand{\im}{{\operatorname{im}}\,}
\newcommand{\Aut}{{\operatorname{Aut}}}

\newcounter{Tc}

\newtheorem{thm}[Tc]{Theorem}
\newtheorem{prop}[Tc]{Proposition}

\newtheorem{cor}[Tc]{Corollary}

\title{Gr\"un's Theorems and Class Groups}
\author{Franz Lemmermeyer}
\address{M\"orikeweg 1, 73489 Jagstzell, Germany}
\email{hb3@ix.urz.uni-heidelberg.de}
\subjclass{11R29}
\keywords{ideal class group; Galois group}
\begin{document}

\begin{abstract}
In this article we show how Gr\"un's results in group theory
can be used for studying the structure of class groups in normal 
extensions.
\end{abstract}

\maketitle

Otto Gr\"un, who is known for his contributions to group theory,
was an amateur mathematician; for details on his ``remarkable career'' 
and his correspondence with Hasse, see Roquette's article \cite{Roq}.
Gr\"un's idea of studying the action of the Galois group on class 
groups led him into problems in representation theory and group theory. 
In this article, I would like to explain how Gr\"un's results
can be applied to class field theory.

\section{Gr\"un's Problem 153}

Gr\"un's excursions into group theory started with the 
following observation\footnote{Gr\"un used the rationals 
as his base field, but everything is valid over an arbitrary number 
field.}, which Gr\"un communicated to Hasse in a letter dated 
Dec. 19, 1932:

\begin{prop}\label{PG0}
Let $L/K$ be a normal extension with degree $n$ and Galois 
group $G$, and assume that the class group $\Cl(L)$ of $L$ 
is cyclic. If $G$ equals its commutator subgroup $G'$, then 
$c^n = 1$ for all $c \in \Cl(L)$.
\end{prop}

Gr\"un did not trust his results and asked Hasse whether he could 
find an error. Hasse could not, but suggested to ask Olga Taussky 
whether the result was known. Taussky sent Arnold Scholz a 
counterexample, which Scholz showed not to be valid; in addition
he proved a large part of Gr\"un's claim by observing that the 
action of a group $G$ on a cyclic group is abelian. Eventually 
Hasse suggested posing Gr\"un's result as a problem in the 
Jahresberichte der DMV. By then, Holzer had given his result 
the following form:

\begin{prop}\label{PG1}
Let $L/k$ be a normal extension of number fields with Galois 
group $G$, and let $K/k$ be the maximal abelian subextension 
of $L/k$. If $\Cl(L)$ is cyclic, then $h_L \mid (L:K)h_K$.
\end{prop}

\section{Galois Action}

Results such as Gr\"un's should be seen as variations of a now
classical theme: the action of a group $G$ on a finite abelian
group $A$ puts contstraints on the structure of $A$. In fact,
define the relative class group $\Cl(K/k)$ of a normal extension
$K/k$ of number fields by the exact sequence
$$ \begin{CD}
   1 @>>> \Cl(K/k) @>>> \Cl(K) @>N>> \Cl(k) 
   \end{CD} $$ 
induced by the norm map $N: \Cl(K) \lra \Cl(k)$.
Then we have

\begin{prop}\label{PGal}
Let $K/k$ be a cyclic extension of prime degree $n$.  For a prime 
$p \nmid n$, let $f$ denote the order of $p \bmod n$. Then the 
$p$-rank of $\Cl(K/k)$ is divisible by $f$.
\end{prop}

The first special case of this proposition was found by Kummer,
who used a similar technique for showing that the class group of
$\Q(\zeta_{29})$ has type $(2,2,2)$. For proofs of more general 
results in this direction see \cite{LGal}; the essential lines in the
historical development were described by Mets\"ankyl\"a \cite{Mets}.

\begin{cor}
If $\Cl_p(K)$ has rank $< f$, then it comes from $k$.
\end{cor}

We say that a piece $C$ of the class group $\Cl(K)$ comes from a 
subfield $k$ if the norm map $N: C \lra \Cl(k)$ is injective.

\begin{proof}
Prop. \ref{PGal} implies that $\Cl_p(K/k)$ is trivial, hence the 
norm map $\Cl_p(K) \lra \Cl_p(k)$ is injective.
\end{proof}

Gr\"un's first attempt at generalizing his results on cyclic class
groups was the following theorem, which he announced in a letter 
to Hasse dated Dec. 5, 1933:

\begin{thm}\label{GT1}
Let $G$ be a finite group acting on an elementary-abelian $p$-group $A$,
and let $m$ denote the rank of $A$. For a prime $\ell \ne p$ let
$m_\ell$ denote the order of $p \bmod \ell$, and let $P$ denote
an $\ell$-Sylow subgroup of $G$.

If $\ell > \frac{m}{m_\ell}$, then $P'$ acts trivially on $A$;
that is, the $\ell$-Sylow subgroups of the automorphism group of $A$
are abelian.
\end{thm}

Gr\"un's Thm. \ref{GT1} does not have nontrivial applications
for $\ell = 2$, since in this case we have $m_2 = 1$, and the
result applies only if $2 > m$, i.e., if $\Cl_2(K)$ is cyclic.
For applications to $p$-class groups for $p \ge 3$, see Section
\ref{SApp} below.
In his article \cite{Gr35}, Gr\"un removed the restriction
that $A$ have exponent $p$, and showed

\begin{thm} \label{GT2}
Let $A$ be an abelian $p$-group on which a finite group $G$ acts. 
If $\nu$ is the smallest exponent with $\ell^\nu > \frac{m}{m_\ell}$, 
then the $\nu$-th commutator subgroup of an $\ell$-Sylow subgroup of 
$G$ acts trivially on $A$.
\end{thm}

Gr\"un's proof is based on the structure of $G_n = \GL_n(\F_q)$, 
which he had obtained in \cite{Gr34}. There he started by observing that 
$$ \# G_n = (p^{nf}-1)(p^{nf}-p^f) \cdots (p^{nf} - p^{(n-1)f}), $$
where $q = p^f$, defined $m_\ell$ as the minimal positive integer
with $p^{fm_\ell} \equiv 1 \bmod \ell$, and then proved

\begin{thm}\label{GT-2}
The $\ell$-Sylow subgroups of $\GL_n(\F_q)$ are elementary 
abelian $\ell$-groups of order $\ell^{ri}$ if
$\big[\frac{n}{m_\ell}\big] = r < \ell$ and 
$\ell^i \parallel (p^{fm_\ell}-1)$.

They are $(r+1)$-stage metabelian groups if 
$$ \ell^r \le \Big[ \frac{n}{m_\ell}\Big] < \ell^{r+1}. $$ 
\end{thm}

Gr\"un next claims (\cite[Satz 1]{Gr35})

\begin{thm}\label{TS1}
Assume that $G$ is a finite group, and that $A \le G$ is an abelian 
$p$-group of $p$-rank $m$ and normal in $G$; let $L$ denote an 
$\ell$-Sylow subgroup of $G$ for some prime $\ell \ne p$, and let 
$m_\ell$ denote the order of $\ell \bmod p$.
 
If $k$ is minimal with $\frac{m}{m_\ell} < \ell^k$, then $A$ commutes 
with the $k$-th commutator group $L^{(k)}$. In other words: $L^{(k)}$ 
acts trivially on $A$.
\end{thm}

For the proof of Theorem \ref{TS1}, Gr\"un considers the map
$\rho: G \lra \Aut(A)$ and assumes first that $A$ is elementary
abelian. Then $\im \rho$ can be identified with some subgroup of 
$\GL_m(p)$. Let $N$ be the group of all elements of $G$ that 
commute with $A$. Then $N$ is normal with 
$G/N \simeq \im \rho \subseteq \GL_m(p)$.  By Thm. \ref{GT-2}, 
the $k$-th commutator group of the $\ell$-Sylow of $G/N$ is
trivial. But the $\ell$-Sylow subgroups of $G/N$ have the form
$LN/N$, where $L$ is an $\ell$-Sylow subgroup of $G$, and we
have $(LN/N)^{(k)} \simeq L^{(k)}N/N$; this shows that 
$L^{(k)} \subseteq N$, which means that $L^{(k)}$ acts
trivially on $A$.
The claim for general abelian $p$-groups $A$ is proved by induction.

\section{Applications}\label{SApp}

In this section we show how Gr\"un's results can be applied to
obtaining results on ideal class groups in normal extensions. 
The general setup is this: $L/k$ is a normal extension with 
Galois group $G$; this group $G$ acts on $\Cl(L)$ or on certain 
subgroups; if some subgroup of $G$ acts trivially, parts of the 
class group must come from proper subfields of $L$.

In fact, let $L/k$ be a nonabelian $\ell$-extension with Galois
group $G$, and let $K$ be the fixed field of $G'$. Let $m$
denote the rank of the $p$-class group of $L$. Gr\"un's Theorem
\ref{GT1} shows:
\begin{enumerate}
\item If $p \equiv 1 \bmod \ell$ and $m = 1$, then $\Cl_p(L)$
      comes from $K$.
\item If $p \equiv -1 \bmod \ell$ ($\ell \ne 2$) and $m \le 2$, 
      then $\Cl_p(L)$ comes from $K$.
\item If $p \not \equiv \pm 1 \bmod \ell$ and $m \le 3$, then  
      $\Cl_p(L)$ comes from $K$.
\end{enumerate}

Consider e.g. an extension $L/k$ whose Galois group $G$ is one
of the two nonabelian groups of order $\ell^3$. Its commutator
subgroup $G'$ is cyclic of order $\ell$, and $G/G' \simeq (\ell,\ell)$.
Let $K$ be the fixed field of $G'$; then $K/k$ is the maximal abelian
subextension of $L/K$, and $\Gal(K/k) \simeq (\ell,\ell)$.

If $A = \Cl_p(L)$ has $p$-rank $m$, and if $m_\ell$ denotes the
order of $p \bmod \ell$, then $G'$ acts trivially if $\ell > m/m_\ell$.
Thus if $\Cl_p(L)$ has $p$-rank at most $2$, then it must come
from $K$. Stronger results follow if there are nontrivial lower bounds 
for $m_\ell$: unless $p \equiv 1 \bmod \ell$, we have $m_\ell \ge 2$, 
and the conclusion above holds as soon as $2\ell > m$.

If $\ell = 3$ and $p \equiv 2 \bmod 3$, then $\ell > m/m_3$ is satisfied 
whenever $m < 6$. Thus if rank $\Cl_2(L) < 6$, then $G = \Gal(L/\Q)$
acts in such a way on the $p$-elementary part of the class group
that the commutator group acts trivially.

\begin{prop}
Let $L/\Q$ be a normal extension whose Galois group is isomorphic to one 
of the two nonabelian groups of order $\ell^3$. If $p \equiv 2 \bmod 3$ 
is prime and if the rank of $\Cl_p(L)$ is $\le 5$, then $\Cl_p(L)$ comes 
from the maximal abelian subfield $K$ of $L/\Q$.
If $p \equiv 1 \bmod 3$ and the rank of $\Cl_p(L)$ is $\le 2$,
then $\Cl_p(L)$ comes from the maximal abelian subfield $K$ of $L/\Q$.
\end{prop}

We have already remarked that Gr\"un's Thm. \ref{GT1} does not apply 
to $2$-extensions. But we can apply Thm. \ref{GT2} to $2$-class field 
towers of sufficiently high stage. Consider e.g. the $2$-class field 
towers of complex quadratic number fields $k$ with class group 
$(2,2,2)$ studied e.g. in \cite{BLS}, \cite{Lem97}, and \cite{Nover}. 
If $k^3$ has a class number divisible by an odd prime $p$, then either 
it comes from $k^2$, or the $p$-rank of $\Cl(k^3)$ is large:

\begin{prop}
Let $k$ be a quadratic number field whose $2$-class field tower
$k \subset k^1 \subset k^2 \subset \ldots$ has at least $3$ steps,
and let $G = \Gal(k^3/k)$ denote the Galois group of the
$2$-extension $k^3/k$. If $C = \Cl_p(k^3)$ has rank $\le 3$, 
then $G''$ acts trivially on $C$, and $\Cl_p(k^2)$ has a subgroup
isomorphic to $C$.
\end{prop}

\end{document}